\documentclass[12pt]{article}

\def\edo{\end{document}	 }

\usepackage{amsmath, amsthm, amscd, amsfonts, amssymb}

\newtheorem{theorem}{Theorem}[section]

\newtheorem{corollary}[theorem]{Corollary}
\newtheorem{lemma}[theorem]{Lemma}

\newtheorem{remark}[theorem]{Remark}

\usepackage[usenames,dvipsnames]{color}

\def\divv{{\rm div }}
\def\rrd{{\mathbb{R}^d}}

\def\calf{{\mathcal{F}}}
\def\calm{{\mathcal{M}}}
\def\calo{{\mathcal{O}}}

\def\cald{{\mathcal{D}}}
\def\calx{{\mathcal{X}}}
\def\calk{{\mathcal{K}}}
\def\call{{\mathcal{L}}}

\def\calp{{\mathcal{P}}}

\def\vsp{\vspace*{1,5mm}\\ }
\def\bk{\bigskip }

\def\sk{\smallskip }
\def\n{\noindent }
\def\dd{\displaystyle}

\def\barr{\begin{array}}
\def\earr{\end{array}}

\def\bit{\begin{itemize}}
\def\eit{\end{itemize}}

\def\FP{Fokker--Planck}

\def\1{^{-1}}

\def\E{{\mathbb{E}}}
\def\rr{{\mathbb{R}}}
\def\nn{{\mathbb{N}}}
\def\9{{\infty}}
\def\lbb{{\lambda}}

\def\ov{\overline}
\def\vf{{\varphi}}
\def\oo{{\omega}}
\def\ooo{{\Omega}}
\def\pp{{\partial}}
\def\vp{{\varepsilon}}

\def\ff{\forall }
\def\({\left(}
\def\){\right)}
\def\<{\left<}
\def\>{\right>}

\title{The ergodicity of nonlinear  Fokker--Planck flows in $L^1(\rrd)$}
\author{Viorel Barbu\thanks{Octav Mayer Institute of Mathematics of  Romanian Academy,  and Al.I. Cuza University,   Ia\c si, Romania.  Email: vbarbu41@gmail.com}\and Michael R\"ockner\thanks{Fakult\"at f\"ur Mathematik, Universit\"at Bielefeld,  D-33501 Bielefeld, Germany.  Email: roeckner@math.uni-bielefeld.de}}
\date{}

\begin{document}
\maketitle
\begin{abstract}
\n One proves in this work that the non\-linear semigroup $S(t)$ in $L^1(\rrd)$, $d\ge 3$, asso\-cia\-ted with the nonlinear Fokker--Planck equation   $u_t-\Delta\beta(u)+\divv(Db(u)u){=}0$, $u(0)=u_0$ in $(0,\9)\times\rrd$,  under suitable conditions on the coefficients $\beta:\rr\to\rr$, $D:\rrd\to\rrd$ and $b:\rr\to\rr$, is mean ergodic. In particular, this implies the mean ergodicity of  the time marginal laws of the solutions to the correspon\-ding McKean--Vlasov stochastic differential equation. This completes the results established in \cite{7} on the nature of the corresponding omega-set $\oo(u_0)$ for $S(t)$  in the case where the   flow $S(t)$ in $L^1(\rrd)$ has not a fixed point and so the corresponding  stationary \FP\ equation has no  solutions.\sk
\end{abstract}

\section{Introduction}\label{s1}
Consider the nonlinear \FP\ equation %(NFPE) 
\begin{equation}\label{e1.1}
\barr{l}
u_t-\Delta\beta(u)+{\rm div}(Db(u)u)=0,\ \mbox{ in  } (0,\9)\times\rr^d,\\
u(0,x)=u_0(x),\ x\in\rr^d,
\earr\end{equation}where $\beta:\rr\to\rr$, $D:\rrd\to\rrd$, $d\ge3$, 
and $b:\rr\to\rr$ are given functions to be made precise in the following.

This equation describes, in statistical physics and mean field theory, the dynamics of a set of interacting particles or of many body systems in disordered media (the so-called anomalous diffusion). (See, e.g., \cite{12}.) In such a situation, for each $t\ge0$, $u=u(t,\cdot)$ is a probability density for each probability density $u_0$.  Another source for equation \eqref{e1.1} is the description of the dynamics of It\^o stochastic processes $X(t)$ in terms of their probability densities $u=u(t,x)$. Namely, if $u\in L^1_{\rm loc}((0,\9)\times\rrd)$ is a Schwartz distributional solution to \eqref{e1.1} such that $t\to u(t,\cdot)dx$ is weakly continuous and $u(t,\cdot)$ is a probability density, then there is a probabilistically weak solution $X$ to the McKean--Vlasov stochastic differential equation in $\rrd$
\begin{equation}\label{e1.2}
	\hspace*{-5mm}\barr{l}
	dX(t)=D(X(t))b
	\(\dd\frac{d\call_{X(t)}}{dx}\,(X(t))\)dt
	+\!\sqrt{\dd\frac{2\beta\(\frac{d\call_{X(t)}}{dx}(X(t))\)}{\dd\frac{d\call_{X(t)}}{dx}(X(t))}}\ dW(t),\\
	\call_{X(t)}(dx)=u(t,x)dx,\ t\ge0, \earr	\hspace*{-5mm}\end{equation}on a probability space $(\ooo,\calf,\mathbb{P})$ with the natural filtration $(\calf_t)_{t\ge0}$, where $W$ is a $d$-dimensional $(\calf_t)$-Brownian motion. Here, $\call_{X(t)}$ is the law of the process $X(t)$ under $\mathbb{P}$ and $\call_{X_0}=u_0dx$ (see \cite{3} for details). 

A function $u:[0,\9)\times\rrd\to\rr$ is called a {\it mild solution} to \eqref{e1.1} if $u\in C([0,\9);L^1(\rrd))$ and
\begin{equation}\label{e1.3}
u(t)=\lim_{h\to0}u_h(t)\mbox{ strongy in }L^1(\rrd),\ \ff t\ge0, \end{equation}where $u_h:[0,\9)\to L^1(\rrd)$ is the solution to the equation
\begin{equation}\label{e1.4}
 \barr{ll}
\dd\frac1h\,(u_h(t)-u_h(t-h))+A_0u_h(t)=0&\mbox{for }t\ge0,\\ u_h(t)=u_0&\mbox{for }t<0,  
	\earr\end{equation}  and $A_0$ is the operator in $L^1(\rrd)$ defined by
\begin{eqnarray}
	A_0(y)&=&-\Delta\beta(y)+\divv(Db(y)y)\mbox{ in }\cald'(\rrd),\ \ff y\in D(A_0),\label{e1.5}\\[1mm]
	D(A_0)&=&\{y\in L^1(\rrd);\ \beta(y)\in L^1_{\rm loc}(\rrd),\ Db(y)y\in L^1_{\rm loc}(\rrd;\rrd);\quad\label{e1.6}\\
	&&\ \ -\Delta\beta(y)+\divv(Db(y)y)\in L^1(\rrd)\}.\nonumber
	\end{eqnarray}
In particular, the mild solution is a distributional solution (in the sense of L.~Schwartz) of equation \eqref{e1.1}. The existence of a mild solution $u$ to \eqref{e1.1} was studied under various hypotheses on $\beta,d$ and $b$ in the works \cite{3}--\cite{7}. The idea, previously used by M.G. Crandall in the existence theory of entropy solutions to a nonlinear conservation law equation \cite{10}, is to represent \eqref{e1.1} as a Cauchy problem in $L^1(\rrd)$
\begin{equation}\label{e1.7}
\frac{du}{dt}+Au=0,\ \ff\,t\ge0;\quad u(0)=u_0,\end{equation}
where $A$ is an $m$-accretive operator in $L^1(\rrd)$ such that $(I+\lbb A)\1f\in(I+\lbb A_0)\1f,$ $\ff f\in ^1(\rrd)$, $\lbb>0.$ Then, by the Crandall \& Liggett generation theorem (see, e.g., \cite{1}, \cite{2}) there exists 
\begin{equation}
	\label{e1.8}
	S(t)u_0=u(t,u_0)=\lim_{n\to\9}\(I+\frac tn\,A\)^{-n}u_0,\ \ \ff t\ge0,\ u_0\in\ov{D(A)},\end{equation}strongly in $L^1(\rrd)$, uniformly in $t$ on bounded intervals. The function $u=u(t,u_0)$ is a mild solution to equation \eqref{e1.1} in the sense of \eqref{e1.3}--\eqref{e1.4} and the mapping $S(t):\ov{D(A)}\to\ov{D(A)}$, $t\ge0$, is a continuous semigroup of contractions in $L^1(\rrd)$ on $\ov{D(A)}$ -- the closure of the domain $D(A)$ of $A$ in $L^1(\rrd)$. We call such a semigroup of contractions a {\it nonlinear \FP\ flow}. It should be emphasized that, in general, this  semigroup $S(t)$ is not unique because its generator $A$ is constructed from $A_0$ by
\begin{equation}\label{e1.9}A(y)=A_0(J_\lbb(y)),\ \ff y\in D(A)=\{u=J_\lbb(f);\ f\in L^1(\rrd),\ \lbb>0\},
	\end{equation}where $J_\lbb:L^1(\rrd)\to L^1(\rrd)$ is a family of contractions such that $J_\lbb(f)\in(I+\lbb A_0)\1f$, $\ff f\in L^1(\rrd)$, $\lbb>0$. Since $I+\lbb A_0$ is, in general, not one-to-one, hence $(I+\lbb A_0)\1$ is multivalued, the family $\{J_\lbb\}_{\lbb>0}$ is not unique, hence so is the operator $A$.  There is an alternative approach to existence theory for nonlinear \FP\ equations deve\-loped by J.A. Carrillo \cite{8}, G.Q. Chen and B. Perthame \cite{9} in the context of entropy and kinetic solutions, but we shall not pursue  this approach in this paper. (In fact, a mild solution to \eqref{e1.1} is a weaker concept of solution than that of entropy solutions.) The semigroups $S(t)$ represents a section in the class of mild solutions.

Here, we shall consider equation \eqref{e1.1} under the following hypotheses:
\begin{itemize}
	\item[\rm(H1)] $\beta\in C^1(\rr),\ \beta'(r)>0,\ \ff r\in\rr\setminus\{0\},\ \beta(0)=0$, and
\begin{equation}\label{e1.9a}
\mu_1\min\{|r|^\nu,|r|\}\le|\beta(r)|\le\mu_2|r|,\ \ff r\in\rr,\end{equation}for $\mu_1,\mu_2>0$ and $\nu>\frac{d-1}d,\ d\ge3$.
	\item[\rm(H2)] $D\in L^\9(\rrd;\rrd)\cap W^{1,1}_{\rm loc}(\rrd;\rrd),\,  \divv\,D\in (L^2(\rrd){+}L^\9(\rrd)),\, D{=}-\nabla\Phi$, where $\Phi\in C(\rrd)\cap W^{1,1}_{\rm loc}(\rrd)$ and 
	$$\barr{c}
	\Phi(x)\ge1,\ \ff x\in\rrd,\ \dd\lim_{|x|\to\9}\Phi(x)=+\9,\vsp 
\Phi^{-m}\in L^1(\rrd\mbox{\ \ for some }m\ge2,\vsp 
\mu_2\Delta\Phi(x)-b_0|\nabla\Phi(x)|^2\le0,\mbox{ a.e. }x\in\rrd.\earr$$
\item[\rm(H3)] $b\in C^1(\rr)\cap C_b(\rr),\ b(r)\ge b_0>0$ for all $r\in[0,\9).$\end{itemize}
We note that hypothesis (H1) does not preclude the degeneracy of the nonlinear diffusion function $\beta$. For instance, any continuous, increasing function $\beta:\rr\to\rr$ of the form
$$\beta(r)=\left\{\barr{ll}
\mu_1r|r|^{d-1}&\mbox{ for } |r|\le r_0,\\ 
\mu_2h(r)&\mbox{ for }|r|>r_0,\earr\right.$$where $r_0>0$, $\mu_1,\mu_2>0,\ |h(r)|\le L|r|,\ \ff r\in\rr,\ L>0$, satisfies \eqref{e1.9a}. As regards Hypothesis (H2),    an example of such a function $\Phi$ is the following
$$\Phi(x)=\left\{\barr{ll}
|x|^2\log|x|+\mu&\mbox{for }|x|\le\delta=\exp\(-\frac{d+2}{2d}\),\\ 
\vf(|x|)+\eta|x|+\mu&\mbox{ for }|x|>\delta,\earr\right.$$where $\mu,\eta>0$ are sufficiently large and $\vf$ is as in \cite[Appendix]{5}. 
As a matter of fact, in this case (see \cite{4}--\cite{7}) the family $\{J_\lbb\}_{\lbb>0}$ of resolvents which defines the operator $A$ is given by the viscosity approximation scheme
\begin{equation}\label{e1.10}
J_\lbb(f)=\lim_{\vp\to0}y_\vp\mbox{\ \ in $L^1(\rrd)$},\end{equation}
where $y_\vp$ is the solution to the equation
\begin{equation}\label{e1.11}
y_\vp-\lbb\Delta(\beta_\vp(y_\vp)+\vp y_\vp)+\lbb\,\divv(D_\vp b_\vp(y_\vp)y_\vp)=f\mbox{\ \ in }\rrd,\end{equation}while  $\beta_\vp,D_\vp$ and $b_\vp$ are smooth approximations of $\beta,D$ and $b$.

However, it turns out  (see \cite{6}) that if one further assumes that, for some $\alpha>0$,
\begin{equation}\label{e1.12}
|b(r)r-b(\bar r)\bar r|\le\alpha|\beta(r)-\beta(\bar r)|,\ \ff r,\bar r\in\rr,\end{equation} then $(I+\lbb A_0)\1$ is single-valued and so $A$ is uniquely defined, more precisely $A=A_0$. In the following, we shall consider the continuous semigroup $S(t)$ gene\-ra\-ted by $A$, which is given by \eqref{e1.9},  and we shall call it the nonlinear \FP\ {\it flow}. This semigroup leaves invariant the set $\calp$ of all the probability densities $\rho$ on $\rrd$, that is,\newpage
$$\calp=\left\{\rho\in L^1(\rrd);\ \rho\ge0,\mbox{ a.e. in }\rrd;\ \int_\rrd\rho\,dx=1\right\}.$$
Now, consider the orbit $\gamma(u_0)=\{S(t)u_0,\ t\ge0\}$ of $S(t)$,  where $u_0\in C:=\ov{D(A)}=L^1(\rrd)$ if $\beta\in C^2(\rr)$ and we associate to $u_0$ the $\oo$-limit~set
\begin{equation}\label{e1.13a}
	\barr{lcl}
	\oo(u_0)&=&\dd
	\left\{u_\9=\lim S(t_n)u_0
	\mbox{ in $L^1$ for some $\{t_n\}\to\infty$}\right\}\vsp
	&=&\dd\bigcap_{s\ge0}\ \  \ov{\bigcup_{t\ge s}S(t)u_0}.\earr
\end{equation}
This set is an attractor for $S(t)$ and, in particular, if $\oo(u_0)\ne\emptyset$ and   it consists of one element $u_\9$ only, then we have
$$\lim_{t\to\9}S(t)u_0=u_\9\mbox{ in }L^1.$$
In \cite{5}, it was proved that, if $\beta$  is not degenerate in the origin, that~is,
\begin{equation}\label{e1.13}
	0<\gamma_0\le\beta'(r)\le\gamma_1,\ \ \ff\,r\in\rr,\end{equation}(which also  implies that $C=L^1$), then, for each $u_0\in \calp$, such that
\begin{equation}\label{e1.14}
	u_0 \ln(u_0)\in L^1(\rr^d),\ \|u_0\|=\int_{\rrd}u_0(x)\Phi(x)dx<\9,
\end{equation}
one has $\oo(u_0)=\{u_\9\}$, where $u_\9$ is an equilibrium state of the flow $S(t)$ and, as a matter of fact, it is the unique solution in $(L^1\cap L^\9)(\rrd)$ to the stationary \FP\ equation
\begin{equation} \label{e1.15a}
	-\Delta\beta(u)+\divv(Db(u)u)=0\mbox{\ \ in }\cald'(\rrd).\end{equation}
 This is an $H$-theorem type result for the \FP\ equation \eqref{e1.1} (see, e.g., \cite{12}, \cite{13} for physical significance and examples). In \cite{7}, the non\-de\-ge\-neracy condition was relaxed to  (H1), which along with (H2)--(H3) leads to the conclusion that, if $u_0\in\calp\cap C$ and $\|u_0\|\le\eta$ for some $\eta>0$, then $\oo(u_0)$ is nonempty, invariant under $S(t)$, $t\ge0$, and compact in $L^1(\rrd)$ which implies that it is an attractor for the trajectory $\gamma(u_0)$. Such a situation occurs if the stationary  equation \eqref{e1.15a} has multiple solutions. If, in addition, there is a fixed point  $a$ of $S(t)$, that is, there is $a\in\calp\cap C$ such that $\|a\|\le\eta$ and $S(t)a=a$ for $t>0$, then $\oo(u_0)$ lies on a ball $\{y\in L^1(\rrd);\ |y-a|_{L^1(\rrd)}=r\}$. In \cite{7}, sufficient conditions on $\beta$ and $b$ for the existence of such a fixed point for $S(t)$ were given. For instance, this happens if
 $$\lim_{r\to+\9}\int^r_1\frac{\beta'(s)}{sb(s)}\ ds=+\9\mbox{ if }\nu\in\mbox{$\(1-\frac1d,1\right]$}$$and
 $$\lim_{r\to0}\int^r_1\frac{\beta'(s)}{sb(s)}\ ds=-\9\mbox{ if }\nu>1.$$Here, no such   condition will be imposed and so it is not clear whether the semigroup $S(t)$ has a fixed point $a$ in $\calp\cap C$, but the nature of the omega-limit set $\oo(u_0)$ will be made clear from the asymptotic properties of the semigroup $S(t)$.  
 Namely, we shall prove that, under the above hypotheses, the flow $t\to S(t)$ is ergodic in $L^1(\rrd)$ (Theorem \ref{t2.1}) and, as a consequence (see Corollary \ref{c2.3}),  the time marginal laws of the probabilistically weak solution $X$ to the McKean--Vlasov SDEs \eqref{e1.2} are mean ergodic, which by our knowledge is a new result in the theory of McKean--Vlasov equations. The proof of Theorem \ref{t2.1} relies on the property of the flow $S(t)$ to be a semigroup of nonlinear contractions in $L^1(\rrd)$ and on  the existence of a unique Haar measure on compact set $\oo(u_0)$.

\bk\noindent{\bf Notation.} $L^p(\rrd)=L^p,\ 1\le p\le\9$, is   the  space of real-valued Lebesgue measurable, $p$-integrable functions on $\rr^d$ with the norm $|\cdot|_p$. The space $L^p(\rrd;\rrd)$ is ana\-lo\-gously defined and $W^{1,1}(\rrd;\rrd)$ is the Sobolev space $\{u\in L^1(\rrd;\rrd);$ $ D_iu_j\in L^1(\rrd),$ $ i=1,...,d; u=(u_j)^d_{j=1}\}$. By $W^{1,1}_{\rm loc}(\rrd;\rrd)$ we denote the corres\-ponding local space. Let $C_b(\rr)$ denote the space of continuous and bounded functions on $\rr$ and  $C^1(\rr)$ the space of continuously differentiable functions on $\rr$.

We recall that an operator $A:\calx\to\calx$, where $\calx$ is a Banach space, is called $m$-accretive if $R(I+\lbb A)=\calx$, $\ff\lbb>0$, and
$$\|u_1-u_2+\lbb(Au_1-Au_2)\|_\calx\ge\|u_1-u_2\|_\calx,\ \ff\lbb>0,\ u_1,u_2\in D(A),$$where $D(A)$ is the domain of $A$ and $R(I+\lbb A)$ is the range of $I+\lbb A$. 
(See, e.g., \cite{1}, \cite{2}.) For each $\eta>0$, we consider the set
$$\calm_\eta:=\left\{u\in L^1;\|u\|=\dd\int_\rrd|u(x)|\Phi(x)dx\le\eta\right\},$$ where $\Phi$ is the potential of $D$ defined as in hypothesis (H2).

\section{The main result}\label{s2}
\setcounter{equation}{0}

Let $S(t):C\to C$, $C=\ov{D(A)}$, be the semigroup generated by the operator $A$ given above by \eqref{e1.9} and, for a given $\eta>0$, let the set
$$\calk:=\calm_\eta\cap C\cap \calp.$$ Everywhere in the following, we shall assume that hypotheses (H1)--(H3) hold. Theorem \ref{t2.1} which follows is the main result. 

\begin{theorem}\label{t2.1}
	Let $\calx$ be a real Banach space and let $F:\calk\to\calx$ be a uniformly continuous mapping. Then, for each $u_0\in\calk$, the set $\oo(u_0)$ is compact in $L^1(\rrd)$ and
	
	\begin{equation}\label{e2.1}
	\lim_{T\to\9}\frac1T\int^T_0F(S(t)u_0)dt=\int_{\oo(u_0)}F(\xi)d\xi,\end{equation}
where $\oo(u_0)$ is endowed with its natural commutative group structure $($recalled below in the proof of the theorem$)$ and	 $d\xi$ is the normalized Haar measure on $\oo(u_0)$. 
	\end{theorem}

The right hand side  of \eqref{e2.1} is the integral of $F$ with respect to the measure $d\xi$ on the topological group $\oo(u_0)$ (see \cite{14a}).

We recall that the Haar measure $\mu$ on a locally compact topological commutative group $G$ is a nonzero Borel measure $\mu$ which is invariant on $G$, that is, $\mu(gS)=\mu(Sg)=\mu(S)$ for any Borel subset $S\subset G$.

A simple example covered by Theorem \ref{t2.1} is $\calx=\rr$ and $F:L^1\to\rr$ is defined by
\begin{equation}\label{e2.2}
	F(u)=\int_\rrd g(x)u(x)dx,\ \ \ff u\in L^1(\rrd),\end{equation}
where $g\in L^\9(\rrd)$. Then, by Theorem \ref{t2.1}, we obviously have

\begin{corollary}\label{c2.2} Under hypotheses {\rm(H1)--(H3)}, for each $u_0\in\calk$ and $g\in L^\9(\rrd)$, we have
\begin{equation}\label{e2.3}
\lim_{T\to\9}\frac1T\int^T_0 dt\int_\rrd g(x)(S(t)u_0)(x)dx=
\int_{\oo(u_0)}\int_\rrd g(x) \xi(x)\,dx\,d\xi.\end{equation}
Furthermore, the semigroup $S(t)$ is mean-ergodic, that is,
\begin{equation}\label{e2.4}
\lim_{T\to\9}\frac1T\int^T_0S(t)u_0dt=\int_{\oo(u_0)}\xi\,d\xi\mbox{ strongly in }L^1(\rrd),\end{equation}
where $d\xi$ is, as above, the normalized Haar measure on $\oo(u_0)$.
\end{corollary}

To get the latter equation, we apply Theorem \ref{t2.1} with $\calx:=L^1(\rrd)$ and $F=\mbox{inclusion map}$. 

 By Corollary \ref{c2.2} it follows in particular  that, under hypotheses (H1)--(H3)   for the nonlinear \FP\ flow $t\to S(t)u_0$ with  $u_0\in\calk$   the classical {\it Boltzmann hypothesis} (see, e.g., \cite{15}, p.~389) is satisfied with   time average $\int_{\oo(u_0)}\xi\,d\xi$ which is the mean of the Haar measure $d\xi$, on $\oo(u_0)$.  (Such a result is related to the Birkhoff ergodic  theorem \cite{7a}, \cite{14b}.)

Now, coming back to the McKean--Vlasov equation \eqref{e1.2}, we get by Corollary \ref{c2.2} the following ergodicity result for the solutions $X(t)$ to \eqref{e1.2}.
\begin{corollary}\label{c2.3} Let $u_0\in\calk$. Then, under hypotheses {\rm(H1)--(H3)} there is a probabilistically weak solution $X$ to \eqref{e1.2}, where $\call_{X_0}=u_0dx$, such that
\begin{equation}\label{e2.5}
		\lim_{T\to\9}\frac1T\int^T_0 \E[g(X(t))]dt=\int_{\oo(u_0)}\int_\rrd g(x)\xi(x)dx\,d\xi,\ \ff g\in L^\9,\end{equation}and, in particular, 
		\begin{equation}
		\label{e2.6}
		\lim_{T\to\9}\frac1T\int^T_0\call_{X(t)}(B)dt=\int_B\(\int_{\oo(u_0)}\xi\,d\xi\)(x)dx,
		\end{equation}for any Borelian set $B\subset\rrd$.
%Furthermore, we have
%	\begin{equation*}\label{e2.6}
	%	\frac1T\int^T_0\call_{X(t)}(dx)dt\to\int_{\oo(u_0)}\xi(x)dx\,d\xi 
	%	\end{equation*}in the weak topology on $\calp$ a.s. $T\to\9$.
\end{corollary}
\begin{remark}\label{r2.4}\rm Due to the degeneracy of the diffusion coefficienty $\beta$, the case $d=2$ is   singular  for the semigroup approach of equation \eqref{e1.1}, namely for the existence of an $m$-accretive realization $A$ of the operator $A_0$ and this is the principal motivation to avoid it here (see~ \cite{7}). However, the case $d=1$ could be treated in a similar way following \cite{4}, but we omit the details.\end{remark}
    
\section{Proofs}\label{s3}
\setcounter{equation}{0}

\n{\bf Proof of Theorem \ref{t2.1}.} The main step of the proof is to show that for each $u_0\in\calk$ the   set $\oo(u_0)$ is compact in $L^1(\rrd)$. For this, it suffices to prove that the orbit $\gamma(u_0)=\{S(t)u_0;\,t\ge0\}$ of the semigroup $S(t)$ is, for $u_0\in\calk$, precompact in the space $L^1$. To this end, we shall mention  first
 the following lemma (see \cite{7}). 
\begin{lemma}\label{l3.1} Let $\eta>0$ arbitrary but fixed. We have
	\begin{eqnarray}
	&\|(I+\lbb A)\1y\|\le\|y\|,\ \ff \lbb>0,\ y\in\calk,\label{e3.1}\\
	&(I+\lbb A)\1(\calk)\subset\calk\cap D(A),\ \ff\lbb>0,\label{e3.2}\\
	&\|S(t)y\|\le\|y\|,\ \ff y\in C,\ t\ge0,\label{e3.3}\\
	&S(t)(\calk)\subset\calk,\ \ff t\ge0.\label{e3.4}\end{eqnarray}\end{lemma}

\n{\bf Proof.} Recalling \eqref{e1.9} and \eqref{e1.10} for \eqref{e3.1}, it suffices to show that
\begin{equation}\label{e3.5a}
\|y_\vp\|\le\|f\|,\ \ff\lbb>0,\ \vp>0,\end{equation}where $y_\vp$ is the solution to equation \eqref{e1.11}. As regards \eqref{e3.5a}, it follows by (H2), via Lemma 3.2 in \cite{5}. By \eqref{e3.1}, it follows also \eqref{e3.2}, while by the exponential formula \eqref{e1.8} one gets \eqref{e3.3} and \eqref{e3.4}.

Now, we consider the restriction $A^*$ of the operator $A$ to $\calk$, that is, the operator
\begin{equation}\label{e3.5}
A^*(y)=A(y),\ \ff y\in D(A^*)=D(A)\cap\calk.\end{equation}It is easily seen that
$\ov{D(A^*)}\subset\calk\subset(I+\lbb A^*)(D(A^*))=R(I+\lbb A^*),$ $\ff\lbb>0,$ and that 
$(I+\lbb A^*)\1=(I+\lbb A)\1$ on $R(I+\lbb A^*)$. By \eqref{e3.5}, it follows also that the operator $A^*$ with the domain $D(A^*)$ is accretive in $L^1$ and, therefore, by \eqref{e1.8} we have
$$S(t)u_0=\lim_{n\to\9}\(I+\frac tn\,A^*\)^{-n}u_0\mbox{\ \ in }L^1,$$for all $u_0\in\ov{D(A^*)}$ (the closure of the domain $D(A^*)$) and uniformly in $t$ on bounded intervals of $(0,\9)$. In other words, $S(t)$ is a continuous semigroup of contractions on $\ov{D(A^*)}$ generated by $A^*$. Then, 
 to show that the trajectory  $\gamma(u_0)$ is precompact in $L^1$, it suffices to check by Theorem 3.1 in \cite{11b} that 
 
 \begin{lemma}\label{l3.2} The operator $(I+\lbb A^*)\1$ is compact on $\calk$ for $\lbb\in(0,\lbb_0)$ where some   $0<\lbb_0<\9$ are small enough.\end{lemma}

 \n{\bf Proof.} This lemma is just Lemma 4.2 in \cite{7}, but since the proof given there was outlined only and this lemma is a crucial step in our proof, we shall prove it here in all details. To this end, we consider  a sequence $\{f_n\}\subset\calk$ such that $\sup_{n\in\nn}|f_n|_1<\9$ and~set
 $$y_n=(I+\lbb A^*_n)\1 f_n=(I+\lbb_nA_n)\1 f_n.$$We have, therefore (see \eqref{e1.10}--\eqref{e1.11}), $ y_n=\dd\lim_{\vp\to0}y^n_\vp\mbox{ in }L^1,$ where $y^\vp_n$ is the solution to the equation
 \begin{equation}
 \label{e3.6}
 y^n_\vp-\lbb\Delta(\beta(y^n_\vp)+\vp y^n_\vp)+\lbb\,{\rm div}(D_\vp b_\vp(y^n_\vp)y^n_\vp)=f_n\mbox{ in }\rrd.
 \end{equation}Letting $\vp\to0$, we have therefore
 \begin{eqnarray}
&y_n-\lbb\Delta\beta(y_n)+\lbb\,{\rm div}(Db(y_n)y_n)=f_n\mbox{ in }\cald'(\rrd),\label{e3.7}\\
&|y_n|_1\le|f_n|_1,\ \ff n\in\nn.\label{e3.8}\end{eqnarray} 	
To get the compactness of the set $\{y_n\}$ in $L^1$, we need some apriori estimates in the Sobolev space $W^{1,p}$. Namely, we shall prove first that $\beta(y_n)\in W^{1,q}_{\rm loc}(\rrd)$ and, for all $R>0$,
\begin{equation}
\label{e3.9}\|\beta(y_n)\|_{L^q(B_R)}+\|\nabla\beta(y_n)\|_{L^q(B_R)}\le C_R(1+|f_n|_1),\ \ff n\in\nn,
\end{equation}
where $q\in\left[1,\frac d{d-1}\right)$. Here, $B_R=\{x;|x|_d<R\}$. To prove this, we shall use some argument from Lemma 2.4 in \cite{6}. Namely, we have, for all $\vf\in C^\9_0(\rrd)$,
\begin{equation}
\label{e3.10}
\Delta(\vf\beta(y_n)=f_1+{\rm div}\,f_2\mbox{\ \ in }\cald'(\rrd),
\end{equation}where
\begin{equation}
\label{e3.11}
\barr{ll}
f_1\!\!\!&=\dd\frac1\lbb\,(y_n-f)\vf-\beta(y_n)\Delta\vf-(D\cdot\nabla\vf)b^*(y_n),\vsp
f_2\!\!\!&=2\beta(y_n)\nabla\vf+D\vf b^*(y_n).\earr
\end{equation}We set $u=\vf\beta(y_n),\ u_\vp\psi_\vp,\ f^\vp_i=f_i*\psi_\vp,$ $i=1,2,$ where $\psi_\vp$ is a standard mollifier, that is,
$$\psi_\vp(x)=\frac1{\vp^d}\,\psi\(\frac x\vp\),\ \psi\in C^\9_0(\rrd),{\rm support}\,\psi\subset\{x;\,|x|\le1\},\int_\rrd\psi(x)dx=1.$$
Let $\calo,\calo'$ be open balls in $\rrd$ centered at zero such that $\ov{\calo'}\subset\calo$ and choose $\vf\in C^2_0(\rrd)$ such that $\vf=1$ on $\calo'$ and $({\rm supp}\,\vf)_\vp\subset\calo,$ $\vp\in(0,1]$, where $({\rm supp}\,\vf)_\vp$ denotes the closed $\vp$-neighbourhood of supp$\,\vf$. By \eqref{e3.10}, we have
$$\Delta u_\vp=f^\vp_1+{\rm div}\,f^\vp_2\mbox{ in }\calo,\ \ u_\vp\in C\9_0(\calo),$$and $u_\vp=u^1_\vp+u^2_\vp$, where $u^1_\vp,u^2_\vp\in C^\9(\calo)\cap C(\ov{\calo})$ are the solutions to 
\begin{eqnarray}
\Delta u^1_\vp=f^\vp_1\mbox{ in }\calo,&&u^1_\vp=0\mbox{ in }\pp\calo,
\label{e3.12}\\
\Delta u^2_\vp={\rm div}\,f^\vp_2\mbox{ in }\calo,&&u^2_\vp=0\mbox{ on }\pp\calo.\label{e3.13}
\end{eqnarray}Then, by standard elliptic estimates, we have
\begin{equation}
\label{3.8}
\|u^1_\vp\|_{W^{1,q}_0(\calo)}\le C\|f^\vp_1\|_{L^1(\calo)}\le C(|y_n|_1+|f_n|_1),\ \ff\vp>0,
\end{equation}where $1\le q<\frac d{d-1}$ and so, by the Sobolev--Galiardo--Nirenberg theorem it follows by \eqref{3.8} that we have
$$\|u^1_\vp\|_p\le C(|f_n|_1+|y_n|_1),\ \ff p\in\mbox{$\left[1,\frac d{d-2}\)$}\mbox{ if }d>2.$$

If $\psi\in L^m(\calo)$, $m>d$, and $\theta\in W^{2,m}(\calo)\cap W^{1,m}_0(\calo)$ is the solution to the Dirichlet problem
$$-\Delta\theta=\psi\mbox{ in }\calo;\ \ \theta=0\mbox{ on }\pp\calo,$$ we have that 
$$\int_\calo u^2_\vp\Delta\theta dx=-\int_\calo f^\vp_2\cdot\nabla dx\le|f^\vp_2|_1\|\nabla\theta\|_\9\le C(|f|_1+|y_n|_1)\|\psi\|_m.$$Taking into account that $|y_n|_1\le|y|_1,$ $\ff n\in\nn,$ this yields
$$\left|\int_\calo u^2_\vp\psi dx\right|\le C(|f|_1+|y|_1)\|\psi\|_m,\ \ff\psi\in L^m(\calo).$$Then, if $\frac1{m'}=1-\frac1m$, it follows by duality that $u^2_\vp\in L^m(\calo)\subset L^q(\calo)$ for all $q\in\left[1,\frac d{d-1}\right)$ and
$$\|u^2_\vp\|_q\le C(|f|_1+|y|_1),\ \ff\vp>0,$$and so, by \eqref{3.8} it follows also that $u^i_\vp\in L^q(\calo),\ i=1,2,$ and
$$\|u^i_\vp\|_q\le C(|f|_1+|y|_1),\ \ff q\in\mbox{$\left[1,\frac d{d-1}\right)$},\ i=1,2.$$Hence,
$$\|u_\vp\|_q\le C(|f_n|_1+|y_n|_1)\le C(|f_n|_1+1),\ \ff\vp>0,\ q\in \mbox{$\left[1,\frac d{d-1}\right)$}.$$Finally, by letting $\vp\to0$ we get
$$\|\vf\beta(y_n)\|_q\le C(|f_n|_1+|y_n|_1)\le C,\ \ff q\in
\mbox{$\left[1,\frac d{d-1}\right)$}.$$(Here and in the following, $C$ is a positive constant independent of $n$.) Because $\vf$ and the corresponding ball $\calo$ are arbitrary, we conclude that $y_n$, $\beta(y_n)\in L^q_{\rm loc}(\rrd)$ and that (for a possible larger $C$, still independent of $\vp$)
$$\|\beta(y_n)\|_{L^q(B_R)}\le C(|f_n|_1+|y_n|_1)\le C,\ \ff q\in \mbox{$\left[1,\frac d{d-1}\right)$}.$$In particular, this implies that
$$\|f_2\|_q\le C(|f_n|_1+|y_n|_1)\le C<\9$$and, therefore,
\begin{equation}
\label{e3.14}
\|f^\vp_2\|_q\le C(f_n|_1+|y_n|_1)\le C<\9,\ q\in \mbox{$\left[1,\frac d{d-1}\right)$},\ n\in\nn.
\end{equation}
Now,we come back to  \eqref{e3.12}--\eqref{e3.13} and note that $\vf\beta(y)=u=u_1+u_2$, where $u_1,u_2$ are solutions to the equations
\begin{eqnarray}
\Delta u_1&=&f_1\mbox{ in }\cald'(\rrd),\label{e3.12'}\\\Delta u_2&=&{\rm div}(f_2)\mbox{ in }\cald'(\rrd),\label{e3.13'}
\end{eqnarray}where $f_1,f_2$ are defined by \eqref{3.8}.

Since $f_1\in L^1$, it follows that $u_1$ can be represented as
$$u_1=-E*f_1\mbox{ in }\rrd,$$where $E(x)\equiv\frac1{(d-2)\oo_d|x|^{d-2}}$ is the fundamental solution to $\Delta$. 

Hence (see, e.g., \cite{12}), $u_1\in M^{\frac d{d-2}}(\rrd)\subset L^p_{\rm loc}(\rrd),\ \ff p\in \mbox{$\left[1,\frac d{d-1}\right)$}$ and $|\nabla u_1|=|\nabla E*f_1|\in M^{\frac d{d-1}}(\rrd)\subset L^p_{\rm loc}(\rrd),$ $\ff p\in\mbox{$\left[1,\frac d{d-1}\right)$}$ with
\begin{equation}
\label{e3.15}
\|\nabla u_1\|_{L^p(B_R)}\le C(|y_1|_1+|f|_1),\ \ff R>0,\ p\in\mbox{$\left[1,\frac d{d-1}\right)$}.
\end{equation}(Here, $M^\ell$ is the Marcinkievicz space of order $\ell$.)

Similarly, we have
$$u_2=-\nabla E*f_2\mbox{\ \ in }\rrd.$$Taking into account that $|\nabla^2E(x)|\le C|x|^{-d}$, $\ff x\ne0$, and that $\nabla^2E$ is homogeneous of order $d$, it follows by the Calderon--Zygmund theorem \cite{9c} and estimate \eqref{e3.11} that
$$|\nabla u^2_\vp|_p\le C|f^\vp_2|_p\le C(|f_n|_1+|y_n|_1)\le C,\ \ff p\in\mbox{$\left[1,\frac d{d-1}\right)$},$$and, after letting $\vp\to0$, this yields
$$\|\nabla u\|_{L^p(B_R)}\le C_R(|f_n|_1+|y_n|_1),\ \ff p\in\mbox{$\left[1,\frac d{d-1}\right)$},$$and so \eqref{e3.9} holds.

Then, by the Kolmogorov compactness theorem it follows that the sequence  $\{\beta(y_n)\}$ is compact in $L^q_{\rm loc}(\rrd)$ and, therefore, on a subsequence $\beta(y_n)\to\eta$, a.e. on $\rrd$ as $n\to\9$. Since $\beta\in C^1$ and $\beta'>0$, it follows that $y_n\to\beta\1(\eta)$, a.e. on $\rrd$. Since $\nu>\frac{d-1}d$, we may choose $q$ close to $\frac d{d-1}$ such that $\nu q\ge1$ and so, by \eqref{e1.10},
$$\mu^q_1\min\{|r|^q,|r|^{1-q}\}\le|\beta(f)|^q,\ \ff r\in\rr.$$This implies that $y_n\to\beta\1(\eta)$ in every $L^1(K)$, $\subset\rrd$. Hence, the set $\{y_n\}$ is compact in $L^1_{\rm loc}$  and so, by \eqref{e3.1} and taking into account that $\{y_n\}\subset\calm_\eta$, it follows that $\{y_n\}$ is compact in $L^1(\rrd)$, as claimed.~$\Box$

\bk\n{\bf Proof of Theorem \ref{t2.1} (continued)} As mentioned earlier, by Lemma \ref{l3.2} it follows that the
 corresponding  $\oo$-limit set $\oo(u_0)$ is compact in $L^1(\rrd)$. It is  also known from the theory of infinite dimensional dynamical system that, for each $t\ge0$, $\oo(u_0)$ is invariant under $S(t)$ which  is an homeomorphism of $\oo(u_0)$ onto $\oo(u_0)$, that is, $S(t)$ is a group on $\oo(u_0)$. Hence, $\oo(u_0)$ can be endowed with a topological commutative group structure with the product\break $g_1\circ g_2=\lim\limits_{n\to\9}S(t^n_1+t^n_2)u_0$, $g_1,g_2\in\oo(u_0)$, where $g_1=\lim\limits_{n\to\9}S(t^n_1)u_0$,\break $g_2=\lim\limits_{n\to\9}S(t^n_2)u_0$ and $\lim\limits_{n\to\9}t^n_i=+\9$, $i=1,2$. Then, by the classical A.~Weil theorem (see \cite{14a}), there is a unique normalized Haar measure on $\oo(u_0)$ and so, by   Birkhoff's ergodic theorem (see \cite{7a} and  Theorem~1 in \cite{11}) it follows that \eqref{e2.1} holds for each uniformly continuous mapping $F:\calk\to \calx$, and so one obtains Theorem \ref{t2.1}  as a special case.~$\Box$\bk 

\n{\bf Proof of Corollary \ref{c2.3}.}  The existence of a probabilistically weak solution to \eqref{e1.2} follows by \cite[Section 2]{3}, which in turn is based on the superposition principle for  linear \FP\ equations (see \cite{15b}, Theorem 2.5).   Furthermore, formula \eqref{e2.5} follows then by \eqref{e2.3} taking into account~that
$$\E[g(X(t))]=\int_\rrd g(x)(S(t)u_0)(x)dx,\ \ff t\ge0,\ g\in L^\9.$$

 \n{\bf Acknowledgement.} This work was supported by the DFG through SFB 1283/2 2021-317210226 and by a grant of the Ministry of Research, Innovation and Digitization, CNCS--UEFISCDI project  PN-III-P4-PCE-2021-0006, within PNCDI III.


\begin{thebibliography}{nn}
 
\bibitem{1} Barbu, V., {\it Nonlinear Differential Equations of Monotone Type in Banach Spaces}, Springer, Berlin. Heidelberg. New York,  2010.\vspace*{-1,5mm} 

\bibitem{2} Barbu, V., {\it Semigroup Approach to Nonlinear Diffusion Equations}, World Scientific, London, Singapore, Beijing, Hong Kong, Tokyo,  2021.\vspace*{-1,5mm}  

\bibitem{3} Barbu, V., R\"ockner, M.,  From \FP\ equations to solutions of distribution dependent SDE, {\it Annals of Probability}, 48 (2020), 1902-1920.\vspace*{-1mm} 

\bibitem{4} Barbu, V., R\"ockner, M., Solutions for nonlinear \FP\ equations with measures as initial data and McKean-Vlasov equations, {\it J.~Functional Anal.}, 280 (7) (2021), 1-35.\vspace*{-1,5mm}  

\bibitem{5} Barbu, V., R\"ockner, M., {\it The evolution to equilibrium of solutions to nonlinear Fokker-Planck equations}, {\it Indiana University Math. Journal}, 72 (1)  (2023), 89-131.\vspace*{-1,5mm}
 

\bibitem{6} Barbu, V., R\"ockner, M., Uniqueness for nonlinear  Fokker-Planck equations and for McKean-Vlasov SDEs, The degenerate case, {\it J. Funct. Anal.} (to appear).\vspace*{-2mm}   

\bibitem{7} Barbu, V., R\"ockner, M., The invariance principle for nonlinear  Fokker-Planck equations,  {\it J.~Diff. Equations}, 315 (2022), 200-221.\vspace*{-2mm} 

\bibitem{7b} B\'enilan, Ph., Brezis, H., Crandall, M.G., A semilinear elliptic equation in $L^1(\rrd)$, {\it Ann. Scuola Normale Sup. Pisa IV}, vol. II (1975), 523-555.\vspace*{-2mm}

\bibitem{7a} Birkhoff, G., What is the ergodic theorem, {\it Amer. Math. Monthly}, 49 (1942).\vspace*{-2mm} 
 
 \bibitem{9c} Calderon, A.P., Zygmund, A., On singular integrals, {\it American J. Math.}, 78 (1956), 289-300.\vspace*{-1mm}
 	
 \bibitem{8} Carrillo, J., Entropy solutions for nonlinear degenerate problems, {\it Archives Rat. Mech. Anal.}, 147 (1999), 269-361.\vspace*{-2mm} 
 
 \bibitem{9} Chen, G., Perthame, B., Well-posedness for nonisotropic dedgenerate parabolic-hyperbolic equations, {\it Ann. Inst. H. Poincar\'e}, 20 (4) (2003), 645-668.\vspace*{-2mm} 
 
 \bibitem{10} Crandall, M.G., The semigroup approach to first order quasi-linear equations in several space variables,   {\it Israel J. Math.}, 12 (1972), 108-122.\vspace*{-2mm} 
 
 \bibitem{11b}  Dafermos, C., Slemrod, M., Asymptotic behaviour of nonlinear contractions semigroups, {\it J. Fnct. Anjal.}, 13 (1973), 97-110.\vspace*{-2mm}
 
  
 \bibitem{12} Franck, T.D., {\it Nonlinear \FP\ Equations},  Springer, Berlin. Heidelberg. New York, 2005.\vspace*{-2mm}  
 
 \bibitem{13} Franck, T.D., Daffertshofer, A., $H$-theorem for nonlinear \FP\ equations related to generalized thermostatics, {\it Physica A. Statitical Mechanics and Its Applications}, 295 (2001), 455-474.\vspace*{-2mm}   
 
\bibitem{11} Gutman, S., Pazy, A., An ergodic theorem for semigroups of contractions, {\it Proceedings Amer. Math. Soc.}, 88 (1983), 254-256.\vspace*{-2mm} 
 
\bibitem{14b} Nemystskii, V.V., Stepanov, V.V., {\it Qualitative Theory of Differential Equations} (Russian), OGIZ, Moscow, Leningrad, 1947.\vspace*{-2mm}

\bibitem{15b} Trevisan, D., Well-posedness of multidimensional diffusion processes with weakly differentiable coefficients, {\it Electron J. Probab.,} 21 (2016), 1-41.\vspace*{-2mm} 

 
 \bibitem{14a} Weil, A., L'Int\'egration dans des groupes topologiques et ses applications, {\it Actualit\'es Scientifiques et Industrielles}, vol. 869, Hermann, Paris, 1940.\vspace*{-2mm} 
 
 \bibitem{15} Yosida, K., {\it Functional Analysis}, Springer-Verlag, Berlin. Heidelberg. New York, 1971.
 
\end{thebibliography}
\end{document}